\documentclass{article}

\usepackage{mathrsfs,enumerate,amsmath}
\usepackage{amssymb}

 \newtheorem{thm}{Theorem}[section]
 \newtheorem{cor}[thm]{Corollary}
 \newtheorem{lem}[thm]{Lemma}

 \newtheorem{rem}[thm]{Remark}

  \newtheorem{conj}[thm]{Conjecture}
 \numberwithin{equation}{section}

\begin{document}

\title{Special functions associated with positive linear operators}
\author{Ioan Ra\c{s}a}
\maketitle

Subjclass: 33C05, 33C45, 34B40, 41A36

\maketitle

\begin{abstract}
Many well-known positive linear operators (like Bernstein, Baskakov, Sz\'{a}sz-Mirakjan) are constructed by using specific fundamental functions. The sums of the squared fundamental functions have been objects of study in some recent papers. We investigate the relationship between these sums and some special functions. Consequently, we get integral representations and upper bounds for the sums. Moreover, we show that they are solutions to suitable second order differential equations. In particular, we provide polynomial or rational solutions to some Heun equations.
\end{abstract}

\maketitle

\section{Introduction}

The degree of non-multiplicativity for many classical positive linear operators can be described in terms of Chebyshev-Gr\"{u}ss inequalities. Several recent papers have been devoted to this problem; see~\cite{2}, \cite{3}, \cite{7}, \cite{8}, \cite{16} and the references therein. It has been shown in \cite{8} that in studying the degree of non-multiplicativity a significant role is played by the sum of the squared fundamental functions involved in the construction of the operators. Some properties of these sums have been investigated in \cite{8}, \cite{14}, \cite{6}.

In this paper we consider a family of sequences of positive linear operators, depending on a real parameter $c$. The choices $c=-1$, $c=0$, $c=1$ correspond, respectively, to the Bernstein, Sz\'{a}sz-Mirakjan, and Baskakov operators; for details see, e.g., \cite{10}, \cite{5}, \cite{18} and the references therein. The sum $S_{n,c}$ of the corresponding fundamental functions can be expressed in terms of some special functions: the hypergeometric function and the Legendre polynomials. By using this relationship we get integral representations and upper bounds for $S_{n,c}$; several known results of this type are extended and improved. Moreover, we show that $S_{n,c}(x)$ is a solution to a suitable second order differential equation. In particular, for $c\neq 0$, $S_{n,c}\left ( -\frac{x}{c} \right ) $ is a solution of a Heun equation; $S_{n,-1}(x)$ is a polynomial solution, and $S_{n,1}(-x)$ is a rational solution. Concerning Heun equations, their polynomial or rational solutions, and their applications to Physics, see \cite{9}, \cite{11}, \cite{12}, \cite{17}, \cite{19} and the references therein.

\section{Preliminary results}

In this section we recall some results from \cite{5}.

Let $c\in \mathbb{R}$. Set $I_c = \left [ 0, - \frac{1}{c} \right ]$ if $c < 0$, and $I_c = [0, +\infty)$ if $c \geq 0$. For $\alpha \in \mathbb{R}$ and $k \in \mathbb{N}_0$ the binomial coefficients are defined as usual by
\begin{equation*}
{\alpha \choose k}:=\frac{\alpha (\alpha - 1)\dots (\alpha - k+1)}{k!} \quad \mbox{if } k \in \mathbb{N}, \mbox{ and } {\alpha \choose 0}:=1.
\end{equation*}
In particular, ${m\choose k}=0$ if $m\in \mathbb{N}$ and $k > m$.

Let $n>0$ be a real number, $k\in \mathbb{N}_0$ and $x\in I_c$. Define
\begin{equation*}
p_{n,k}^{[c]}(x):=(-1)^k {-\frac{n}{c}\choose k}(cx)^k (1+cx)^{-\frac{n}{c}-k}, \quad \mbox{ if } c\neq 0,
\end{equation*}
\begin{equation*}
p_{n,k}^{[0]}(x):=\lim _{c\to 0}p_{n,k}^{[c]}(x)=\frac{(nx)^k}{k!}e^{-nx}, \quad \mbox{ if } c = 0.
\end{equation*}

Then $\sum _{k=0}^{\infty} p_{n,k}^{[c]}(x) = 1$. Throughout the paper we shall suppose that $n>c$ if $c \geq 0$, or $n=-cl$ with some $l\in \mathbb{N}$ if $c <0$. Under this hypothesis define
\begin{equation}
T_{n,c}(x,y):=\sum _{k=0}^{\infty} p_{n,k}^{[c]}(x)p_{n,k}^{[c]}(y), \quad x, y \in I_c .\label{eq:2.1}
\end{equation}

\begin{lem}
(\cite[Lemma 1]{5}). The series in \eqref{eq:2.1} converges for all $x,y \in I_c$ and
\begin{equation}
T_{n,c}(x,y) = \left ( (1+cx)(1+cy) \right )^{-\frac{n}{c}} {_2F_1}\left ( \frac{n}{c};\frac{n}{c};1;\frac{c^2xy}{(1+cx)(1+cy)} \right ), \quad c\neq 0 ,\label{eq:2.2}
\end{equation}
\begin{equation}
T_{n,0}(x,y) = e^{-n(x+y)} I_0 \left ( 2n\sqrt{xy} \right ), \quad c=0.\label{eq:2.3}
\end{equation}
\end{lem}
Here ${_2F_1}$ is the hypergeometric function and $I_0$ is the modified Bessel function of first kind of order $0$; see \cite[15.1.1 and 9.6.12]{1}.

From the proof of Theorem 1 in \cite{5} we derive:
\begin{thm}
\cite[Theorem 1]{5}. Let $n$ and $c$ be as above, $x,y\in I_c$. For $c \neq 0$ we have
\begin{eqnarray}\label{eq:2.4}
&&T_{n,c}(x,y) = \\ &&\frac{1}{\pi}\int _0^1 \frac{\left ( \left ( \sqrt{c^2xy} + \sqrt{(1+cx)(1+cy)} \right )^2 - 4 t \sqrt{c^2xy(1+cx)(1+cy)} \right )^{-\frac{n}{c}}}{\sqrt{t(1-t)}}dt.\nonumber
\end{eqnarray}
Moreover,
\begin{equation}
T_{n,0}(x,y) = \frac{1}{\pi}\int _{-1}^1 e^{-n (x+y+2t\sqrt{xy})}\frac{dt}{\sqrt{1-t^2}}.\label{eq:2.5}
\end{equation}
\end{thm}

\section{The function $S_{n,c}$}

The functions $p_{n,k}^{[c]}$ are the fundamental functions in the construction of several classes of positive linear operators: see \cite{5}, \cite{10}, \cite{18} and the references therein. In this paper we are interested in the function
\begin{equation}
S_{n,c}(x):=\sum _{k=0}^\infty \left ( p_{n,k}^{[c]}(x) \right )^2, \quad x \in I_c. \label{eq:3.1}
\end{equation}

For particular values of $c$ it was considered in \cite{8} in relation with the problem of studying the non-multiplicativity of some classical positive linear operators. Obviously
\begin{equation}
S_{n,c}(x) = T_{n,c}(x,x), \quad x \in I_c. \label{eq:3.2}
\end{equation}

Consequently, some properties of the function $S_{n,c}$ can be derived from the properties of $T_{n,c}$. In particular, we have:
\begin{thm}\label{thm:3.1}
Let $c\neq 0$. Then
\begin{equation}
S_{n,c}(x) = (1+cx)^{-\frac{2n}{c}} \sum _{k=0}^\infty \left ( \frac{n(n+c)\dots (n+(k-1)c)}{k!} \right )^2 \left ( \frac{x}{1+cx} \right )^{2k}, \label{eq:3.3}
\end{equation}
\begin{equation}
S_{n,c}(x) = (1+cx)^{-\frac{2n}{c}} {_2F_1}\left ( \frac{n}{c};\frac{n}{c};1;\left ( \frac{cx}{1+cx} \right )^2 \right ),\label{eq:3.4}
\end{equation}
\begin{equation}
S_{n,c}(x) = \frac{1}{\pi}\int _0 ^1 \left ( t+(1-t)(1+2cx)^2 \right )^{-\frac{n}{c}} \frac{dt}{\sqrt{t(1-t)}}.\label{eq:3.5}
\end{equation}

If $c=0$, then
\begin{equation}
S_{n,0}(x) = e^{-2nx}\sum _{k=0}^\infty \frac{(nx)^{2k}}{(k!)^2},\label{eq:3.6}
\end{equation}
\begin{equation}
S_{n,0}(x) = e^{-2nx}I_0(2nx),\label{eq:3.7}
\end{equation}
\begin{equation}
S_{n,0}(x) = \frac{1}{\pi} \int _{-1}^{1}e^{-2nx(1+t)} \frac{dt}{\sqrt{1-t^2}}.\label{eq:3.8}
\end{equation}
\end{thm}
\textbf{Proof.}
It suffices to combine \eqref{eq:3.2} with \eqref{eq:2.1}-\eqref{eq:2.5}.
\newline \newline
Let $c\neq 0$. Consider the function
\begin{equation}
H_{n,c}(x):=S_{n,c}\left ( - \frac{x}{c}\right )\label{eq:3.9}
\end{equation}
defined on $(-\infty, 0]$ if $c>0$, and on $[0,1]$ if $c<0$.

The next result shows that $S_{n,c}$ and $H_{n,c}$ satisfy suitable differential equations. Concerning the Heun equations, see \cite{9}, \cite{11}, \cite{12}, \cite{17}, \cite{19} and the references therein.

\begin{thm}\label{thm:3.2}
\begin{enumerate}[(i)]
\item{}The function $S_{n,c}$ is a solution to the differential equation
\begin{eqnarray}
x(1+cx)(1+2cx)y''(x)+(4(n+c)x(1+cx)+1)y'(x)+\label{eq:3.10}\\2n(1+2cx)y(x)=0, \quad x \in I_c. \nonumber
\end{eqnarray}
\item{}If $c\neq 0$, the function $H_{n,c}$ is a solution to the Heun equation
\begin{equation}
y''(x)+\left ( \frac{1}{x} + \frac{1}{x-1}+\frac{\frac{2n}{c}}{x-\frac{1}{2}} \right ) y'(x) + \frac{\frac{2n}{c}x - \frac{n}{c}}{x(x-1)\left ( x-\frac{1}{2} \right )}y(x) = 0, \label{eq:3.11}
\end{equation}
with parameters $\alpha = 1$, $\beta = \frac{2n}{c}$, $\gamma = 1$, $\delta = 1$, $\epsilon = \frac{2n}{c}$, $q=\frac{n}{c}$.
\end{enumerate}
\end{thm}
\textbf{Proof.}
\begin{enumerate}[(i)]
\item{}Let $c\neq 0$ and $n$ be given. The function $w(z):={_2F_1}\left ( \frac{n}{c},\frac{n}{c};1;z \right )$ satisfies (see \cite[15.5.1]{1})
\begin{equation}
c^2z(1-z)w''(z)+c(c-(2n+c)z)w'(z)-n^2w(z)=0. \label{eq:3.12}
\end{equation}

According to \eqref{eq:3.4} we have also
\begin{equation}
w\left ( \left ( \frac{cx}{1+cx}\right )^2\right ) = (1+cx)^{\frac{2n}{c}}S_{n,c}(x). \label{eq:3.13}
\end{equation}

Now it is a matter of calculus to combine \eqref{eq:3.12} and \eqref{eq:3.13} in order to get \eqref{eq:3.10} with $c\neq 0$. Moreover, the function $I_0(z)$ satisfies (see \cite[9.6.1]{1})
\begin{equation}
zI_0''(z) + I_0'(z) - zI_0(z) = 0. \label{eq:3.14}
\end{equation}

According to \eqref{eq:3.7},
\begin{equation}
I_0(2nx) = e^{2nx}S_{n,0}(x). \label{eq:3.15}
\end{equation}

Combining \eqref{eq:3.14} and \eqref{eq:3.15} it is easy to deduce \eqref{eq:3.10} with $c=0$.

\item{}For $c\neq 0$ it is again a matter of calculus to combine \eqref{eq:3.9} and \eqref{eq:3.10} in order to get \eqref{eq:3.11}.
\end{enumerate}

\section{The Bernstein and the Bleimann-Butzer-Hahn bases\label{sect:4}}

The functions $p_{n,k}^{[-1]}(x) = {n \choose  k}x^k (1-x)^{n-k}$, $k=0,1\dots n$, $x\in [0,1]$ are the fundamental Bernstein polynomials, or shortly the Bernstein basis; see \cite[5.2.5]{4}. For the sake of simplicity we shall use the notation
\begin{equation}
F_n (x):=S_{n,-1}(x) = \sum _{k=0}^n \left ( {n \choose k} x^k (1-x)^{n-k}\right )^2, \quad x\in [0,1]. \label{eq:4.1}
\end{equation}

Consider the Legendre polynomials (see \cite[22.3.1]{1}):
\begin{equation}
P_n(x) = 2^{-n}\sum _{k=0}^n {n \choose k}^2 (x+1)^k (x-1)^{n-k}. \label{eq:4.2}
\end{equation}

Let
\begin{equation}
x \in \left [ 0, \frac{1}{2}\right ), \quad t = \frac{2x^2-2x+1}{1-2x} \in [1, +\infty ). \label{eq:4.3}
\end{equation}

Thorsten Neuschel (see \cite{8}) proved that
\begin{equation}
F_n(x) = \left ( t - \sqrt{t^2-1} \right )^n P_n(t), \label{eq:4.4}
\end{equation}
and inferred that $F_n$ is decreasing on $\left [ 0, \frac{1}{2}\right ]$ and increasing on $\left [ \frac{1}{2}, 1\right ]$; so he solved a problem raised in \cite{8}. Another problem was solved by Geno Nikolov; using \eqref{eq:4.4} he proved the following theorem and derived new inequalities involving the Legendre polynomials.

\begin{thm}\label{thm:4.1}
(\cite{14}). $F_n$ is a convex function.
\end{thm}

This result is also a consequence of the following representation:
\begin{equation}
F_n(x) = 4^{-n}{2n \choose n}\sum_{k=0}^n 4^k {n \choose k} ^2 {2n \choose 2k}^{-1} \left ( x-\frac{1}{2} \right )^{2k}. \label{eq:4.5}
\end{equation}

This formula was proved in \cite{6} by using Parseval's formula. Here we present a new proof of \eqref{eq:4.5}; see also Remark \ref{rem:4.4}. To this end, set in \eqref{eq:4.4} $t=\cos \theta$ and $x = \frac{1-e^{-i\theta}}{2}$. Then
\begin{equation}
F_n(x) = \left ( t - \sqrt{t^2-1}\right )^n P_n(t) = e^{-n\theta i}P_n(\cos \theta). \label{eq:4.6}
\end{equation}

On the other hand (see \cite[22.3.13]{1}, \cite[4.9.3]{15}),
\begin{equation}
P_n(\cos \theta) = 4^{-n} {2n \choose n}\sum_{k=0}^n {n \choose k}^2 {2n \choose 2k}^{-1} e^{(n-2k)\theta i}. \label{eq:4.7}
\end{equation}

Now \eqref{eq:4.5} is a consequence of \eqref{eq:4.6} and \eqref{eq:4.7}.

The next result provides another proof (and a slight generalization) of Theorem \ref{thm:4.1}.

\begin{thm}\label{thm:4.2}
For $c<0$, $S_{n,c}$ is a convex function.
\end{thm}
\textbf{Proof.}
For $c<0$ we have $n=-cl$ with some $l \in \mathbb{N}$. According to \eqref{eq:3.5},
\begin{equation}
S_{n,c}(x) = \frac{1}{\pi} \int _0^1 \left ( t + (1-t) (1+2cx)^2 \right )^l \frac{dt}{\sqrt{t(1-t)}}, \quad x\in \left [ 0, -\frac{1}{c} \right ]. \label{eq:4.8}
\end{equation}

It is elementary to verify that $S_{n,c}'' \geq 0$, and this concludes the proof.
\newline \newline

Other properties of the function $F_n$ are presented in the next theorem.

\begin{thm}\label{thm:4.3}
\begin{enumerate}[(a)]
\item{}$F_n$ satisfies the following relations:
\begin{equation}
F_n(x) = \frac{1}{\pi}\int _0 ^1 \left ( t + (1-t)(1-2x)^2 \right )^n \frac{dt}{\sqrt{t(1-t)}}, \label{eq:4.9}
\end{equation}
\begin{equation}
2(n+1)F_{n+1}(x)-(2n+1)(1+(1-2x)^2)F_n(x) + 2n(1-2x)^2F_{n-1}(x) = 0, \label{eq:4.10}
\end{equation}
\begin{eqnarray}\label{eq:4.11}
(1-2x)\left ( F_{n+1}'(x) - (1-2x(1-x))F_n'(x) \right ) = \\ 2(n+2x(1-x)) F_n(x)-2(n+1)F_{n+1}(x) \nonumber,
\end{eqnarray}
\begin{equation}
F_{n+1}' (x) - (1-2x)^2 F_{n-1}'(x) = 2(1-2x)\left ( (2n-1) F_{n-1}(x)-(2n+1)F_n(x) \right ) . \label{eq:4.12}
\end{equation}

\item{}$F_n$ is a solution to the differential equation
\begin{equation}
x(1-x)(1-2x)y''(x) + (1+4(n-1)x(1-x))y'(x)+2n(1-2x)y(x) = 0. \label{eq:4.13}
\end{equation}

Equivalently, $F_n$ is a polynomial solution to the Heun equation
\begin{equation}
y''(x) + \left ( \frac{1}{x} + \frac{1}{x-1} + \frac{-2n}{x-\frac{1}{2}} \right )y'(x) + \frac{-2nx-(-n)}{x(x-1)\left ( x-\frac{1}{2} \right )} y(x) = 0, \label{eq:4.14}
\end{equation}
with parameters $\alpha = 1$, $\beta = -2n$, $\gamma = 1$, $\delta =1$, $\epsilon = -2n$, $q=-n$.

\item{} The following inequality holds:
\begin{equation}
F_n(x) \leq \frac{1}{\sqrt{1+4(n-1)x(1-x)}}, \quad n \geq 1, x\in [0,1]. \label{eq:4.15}
\end{equation}
\end{enumerate}
\end{thm}

\textbf{Proof.}
\begin{enumerate}[(a)]
\item{} \eqref{eq:4.9} is a consequence of \eqref{eq:3.5}. From \eqref{eq:4.3} we deduce that
\begin{equation}
x= \frac{1-t+\sqrt{t^2-1}}{2}, \label{eq:4.16}
\end{equation}
\begin{equation}
t-\sqrt{t^2-1} = 1-2x, \label{eq:4.17}
\end{equation}
\begin{equation}
\frac{dx}{dt} = \frac{(1-2x)^2}{4x(1-x)}. \label{eq:4.18}
\end{equation}

Now \eqref{eq:4.4} and \eqref{eq:4.17} imply
\begin{equation}
P_n(t) = (1-2x)^{-n}F_n(x). \label{eq:4.19}
\end{equation}

On the other hand (see \cite[22.7.1]{1}),
\begin{equation}
(n+1)P_{n+1}(t)-(2n+1)tP_n(t)+nP_{n-1}(t)=0. \label{eq:4.20}
\end{equation}

From \eqref{eq:4.3}, \eqref{eq:4.19} and \eqref{eq:4.20} it is easy to deduce \eqref{eq:4.10}.

As a consequence of \eqref{eq:4.18} and \eqref{eq:4.19} we have
\begin{equation}
P_n'(t) = \frac{(1-2x)^{-n+1}}{4x(1-x)} \left ( (1-2x)F_n'(x) + 2n F_n(x)\right ). \label{eq:4.21}
\end{equation}

The following relations are satisfied by the Legendre polynomials (see \cite[(4.37)]{13}):
\begin{equation}
P_{n+1}'(t) - tP_n'(t) = (n+1)P_n(t), \label{eq:4.22}
\end{equation}
\begin{equation}
P_{n+1}'(t) - P_{n-1}'(t) = (2n+1)P_n(t). \label{eq:4.23}
\end{equation}

Combining \eqref{eq:4.3}, \eqref{eq:4.19}, \eqref{eq:4.21} and \eqref{eq:4.22} we get \eqref{eq:4.11}. To prove \eqref{eq:4.12} it suffices to use \eqref{eq:4.23}, \eqref{eq:4.21}, \eqref{eq:4.19} and \eqref{eq:4.10}.

\item{} \eqref{eq:4.13} and \eqref{eq:4.14} are consequences of \eqref{eq:3.10} and \eqref{eq:3.11}.

\item{} According to Theorem \ref{thm:4.1}, $F_n'' \geq 0$, so that from \eqref{eq:4.13} we infer
\begin{equation*}
\left ( 1+4(n-1)x(1-x) \right ) F_n'(x) + 2n (1-2x) F_n(x)  \leq 0, x\in \left [ 0, \frac{1}{2} \right ].
\end{equation*}

For $n\geq 2$ this leads to
\begin{equation*}
F_n(x) \leq \left ( 1+4(n-1)x(1-x) \right ) ^{-\frac{n}{2(n-1)}}, \quad x \in \left [ 0, \frac{1}{2} \right ].
\end{equation*}

Using the symmetry with respect to $\frac{1}{2}$, we get
\begin{equation}
F_n(x) \leq \left ( 1+4(n-1)x(1-x) \right ) ^{-\frac{n}{2(n-1)}}, \quad x\in [0,1], n\geq 2. \label{4.24}
\end{equation}

This obviously implies \eqref{eq:4.15} for $n\geq 2$. Since $F_1(x)\leq 1$, \eqref{eq:4.15} is valid also for $n=1$, and this concludes the proof.
\end{enumerate}

\begin{rem}\label{rem:4.4}
\eqref{eq:4.5} can be derived from \eqref{eq:4.9} by an elementary calculation.
\end{rem}

The rest of this section is devoted to the Bleimann-Butzer-Hahn basis, which consists of the functions (see \cite[5.2.8]{4})
\begin{equation*}
{n \choose k} x^k (1+x)^{-n}, \quad x \in [0,\infty ), \quad k = 0, 1, \dots , n.
\end{equation*}

Denote
\begin{equation}
U_n(x) = \sum _{k=0}^n \left ( {n \choose k} x^k (1+x)^{-n} \right ) ^2, \quad n \geq 1. \label{eq:4.25}
\end{equation}

It is easy to verify that
\begin{equation}
U_n(x) = F_n \left ( \frac{x}{1+x} \right ), \quad x \in [0, \infty ). \label{eq:4.26}
\end{equation}

Consequently, the properties of the functions $U_n$, presented in the next corollary, can be easily deduced from those of the functions $F_n$.

\begin{cor}
\begin{enumerate}[(a)]
\item{} $U_n$ is a solution to the differential equation
\begin{equation}
x(1-x)(1+x)^2y''(x)+(1+x)(1+4nx-x^2)y'(x)+2n(1-x)y(x)=0. \label{eq:4.27}
\end{equation}
\item{} The following relations hold true:
\begin{equation}
U_n(x) = \frac{1}{\pi} \int _0 ^1 \left ( t+(1-t) \left ( \frac{1-x}{1+x} \right )^2 \right )^n \frac{dt}{\sqrt{t(1-t)}}, \label{eq:4.28}
\end{equation}
\begin{equation}
U_n(x) = 4^{-n} {2n \choose n}\sum _{k=0}^n {n \choose k}^2 {2n \choose 2k}^{-1} \left ( \frac{x-1}{x+1} \right )^{2k}, \label{eq:4.29}
\end{equation}
\begin{equation}
U_n(x)\leq \frac{x+1}{\sqrt{x^2 + (4n-2)x+1}}, \quad x\geq 0, n\geq 1. \label{eq:4.30}
\end{equation}
\end{enumerate}
\end{cor}

\section{The Baskakov basis. The Meyer-K\"{o}nig and Zeller basis}

In this section we are concerned with other two bases, namely (see \cite[5.2.6, (5.3.24)]{4})
\begin{itemize}
\item{}the Baskakov basis:
\begin{equation*}
{n+k-1 \choose k} x^k (1+x)^{-n-k}, \quad x\in [0,+\infty ), k=0,1,\dots ,
\end{equation*}
\item{}the Meyer-K\"{o}nig and Zeller basis:
\begin{equation*}
{n+k \choose k} x^k (1-x)^{n+1}, \quad x\in [0,1), k=0,1,\dots .
\end{equation*}
\end{itemize}

Denote
\begin{equation}
G_n(x) = \sum _{k=0}^\infty \left ( {n+k-1 \choose k} x^k (1+x)^{-n-k} \right )^2, \label{eq:5.1}
\end{equation}
\begin{equation}
J_n(x) = \sum _{k=0}^\infty \left ( {n+k \choose k} x^k (1-x)^{n+1} \right )^2 . \label{eq:5.2}
\end{equation}

It is easy to verify that $G_n = S_{n,1}$ and
\begin{equation}
J_n(x) = G_{n+1} \left ( \frac{x}{1-x} \right ), \quad x\in [0,1), \label{eq:5.3}
\end{equation}
\begin{equation}
G_n(x) = J_{n-1} \left ( \frac{x}{1+x} \right ), \quad x \in [0, +\infty ). \label{eq:5.4}
\end{equation}

It was proved in \cite{6} that
\begin{equation}
J_n(x) = 4^{-n}\sum _{k=0}^n \frac{(2k)!(2n-2k)!}{k!^2(n-k)!^2}\left ( \frac{1-x}{1+x} \right ) ^{2k+1}. \label{eq:5.5}
\end{equation}

It follows that
\begin{equation}
G_n(x) = 4^{1-n}\sum _{k=0}^{n-1} \frac{(2k)!(2n-2k-2)!}{(k!)^2(n-k-1)!} \left ( \frac{1}{2x+1} \right )^{2k+1}. \label{eq:5.6}
\end{equation}

In particular, $G_n$ and $J_n$ are convex functions.

\begin{thm}\label{thm:5.1}
\begin{enumerate}[(a)]
\item{}$G_n$ satisfies the following relations:
\begin{equation}
x(1+x)(1+2x)G_n''(x)+(4(n+1)x(1+x)+1)G_n'(x)+2n(1+2x)G_n(x)=0, \label{eq:5.7}
\end{equation}
\begin{equation}
G_n(x) = \frac{1}{\pi} \int _0 ^1 \left ( t+(1-t)(1+2x)^2 \right )^{-n}\frac{dt}{\sqrt{t(1-t)}}, \label{eq:5.8}
\end{equation}
\begin{equation}
G_n(x) \leq \left ( 4(n+1)x(1+x)+1 \right )^{-\frac{n}{2(n+1)}}. \label{eq:5.9}
\end{equation}
\item{} $G_n(-x)$ is a rational solution to the Heun equation
\begin{equation}
y''(x)+\left ( \frac{1}{x} + \frac{1}{x-1} + \frac{2n}{x-\frac{1}{2}} \right ) y'(x) + \frac{2nx-n}{x(x-1)\left ( x-\frac{1}{2} \right )} y(x) = 0, x<0, \label{eq:5.10}
\end{equation}
with parameters $\alpha = 1$, $\beta = 2n$, $\gamma =1$, $\delta =1$, $\epsilon = 2n$, $q=n$.
\end{enumerate}
\end{thm}
\textbf{Proof.}
Since $G_n = S_{n,1}$, \eqref{eq:5.7}, \eqref{eq:5.8} and \eqref{eq:5.10} follow from Theorems  \ref{thm:3.1} and \ref{thm:3.2}. In order to prove \eqref{eq:5.9}, recall that $G_n$ is convex; combined with \eqref{eq:5.7}, this yields
\begin{equation}
\left ( 4(n+1)x(1+x)+1 \right ) G_n'(x) \leq -2n(1+2x)G_n(x), \quad x\geq 0. \label{eq:5.11}
\end{equation}

Now it is easy to deduce \eqref{eq:5.9} from \eqref{eq:5.11}.

\begin{cor}\label{cor:5.2}
\begin{enumerate}[(a)]
\item{} $J_n$ satisfies the following relations:
\begin{eqnarray}\label{eq:5.12}
x(1+x)(1-x)^2J_n''(x)-(1-x)\left ( x^2-4(n+1)x-1 \right )J_n'(x) +\\ 2(n+1)(1+x)J_n(x) = 0, \nonumber
\end{eqnarray}
\begin{equation}
J_n(x)=\frac{1}{\pi} \int _0^1 \left ( t+(1-t)\left ( \frac{1+x}{1-x} \right )^2 \right )^{-n-1} \frac{dt}{\sqrt{t(1-t)}}, \label{eq:5.13}
\end{equation}
\begin{equation}
J_n(x) \leq \left ( \frac{(1-x)^2}{x^2+(4n+6)x+1} \right )^{\frac{n+1}{2(n+2)}}. \label{eq:5.14}
\end{equation}
\end{enumerate}
\end{cor}

\textbf{Proof.}
It suffices to use \eqref{eq:5.3} and Theorem \ref{thm:5.1}.
\newline \newline

To conclude this section, let us remark that from \cite[(12)]{5} it follows
\begin{equation}
G_n(x) \leq {2n-2 \choose n-1}\frac{(1+x)^{n-1}}{(1+2x)^n}, \quad x\geq 0. \label{eq:5.15}
\end{equation}

Combined with \eqref{eq:5.3}, this yields
\begin{equation}
J_n(x) \leq {2n \choose n} \frac{1-x}{(1+x)^{n+1}}, \quad x \in [0,1). \label{eq:5.16}
\end{equation}

\section{The Sz\'{a}sz-Mirakjan basis}

In this case the fundamental functions are (see \cite[5.3.9]{4})
\begin{equation*}
e^{-nx}\frac{(nx)^k}{k!}, \quad x\in [0, +\infty ), k = 0, 1, \dots .
\end{equation*}

Consequently, we shall consider the function
\begin{equation}
K_n(x):=\sum _{k=0}^\infty \left ( e^{-nx}\frac{(nx)^k}{k!} \right )^2 . \label{eq:6.1}
\end{equation}

\begin{thm}\label{thm:6.1}
The function $K_n$ satisfies:
\begin{equation}
xK_n''(x) + (4nx+1)K_n'(x) + 2nK_n(x) = 0, \label{eq:6.2}
\end{equation}
\begin{equation}
K_n(x) = \frac{1}{\pi}\int _{-1}^1 e^{-2nx(1+t)}\frac{dt}{\sqrt{1-t^2}}, \label{eq:6.3}
\end{equation}
\begin{equation}
K_n(x) \leq \frac{1}{\sqrt{4nx+1}}, \quad x\geq 0. \label{eq:6.4}
\end{equation}
\end{thm}
\textbf{Proof.}
We have $K_n=S_{n,0}$, so that \eqref{eq:6.2} and \eqref{eq:6.3} follow from Theorems \ref{thm:3.1} and \ref{thm:3.2}. From \eqref{eq:6.3} we deduce that $K_n$ is convex (see also \cite{6}) and so \eqref{eq:6.2} implies
\begin{equation*}
(4nx+1)K_n'(x) \leq -2nK_n(x), \quad x\geq 0.
\end{equation*}
This leads immediately to \eqref{eq:6.4} and the proof is finished.

\begin{rem}\label{rem:6.2}
It was proved in \cite{8} that
\begin{equation}
\inf _{x\geq 0} G_n(x) = 0, \quad n\geq 1, \label{eq:6.5}
\end{equation}
\begin{equation}
\inf _{x\geq 0} K_n(x) = 0, \quad n\geq 1. \label{eq:6.6}
\end{equation}

Clearly \eqref{eq:5.9} and \eqref{eq:5.15} are stronger than \eqref{eq:6.5}, and \eqref{eq:6.4} is stronger than \eqref{eq:6.6}.
\end{rem}

\section{Final remarks}

For integral operators of the form $Lf(x) = \int _{I} K(x,y)f(y)dy$, the properties of the function $S(x):= \int _{I}K^2(x,y)dy$ will be presented in a forthcoming paper. There we will discuss similar problems for multivariate operators.
\newline \newline
We conclude this paper with
\begin{conj}
$\log {S_{n,c}}$ is a convex function.
\end{conj}

\end{document}